\catcode`\@=11

\font\tenmsa=msam10
\font\sevenmsa=msam7
\font\fivemsa=msam5
\font\tenmsb=msbm10
\font\sevenmsb=msbm7
\font\fivemsb=msbm5
\newfam\msafam
\newfam\msbfam
\textfont\msafam=\tenmsa  \scriptfont\msafam=\sevenmsa
  \scriptscriptfont\msafam=\fivemsa
\textfont\msbfam=\tenmsb  \scriptfont\msbfam=\sevenmsb
  \scriptscriptfont\msbfam=\fivemsb

\def\hexnumber@#1{\ifnum#1<10 \number#1\else
 \ifnum#1=10 A\else\ifnum#1=11 B\else\ifnum#1=12 C\else
 \ifnum#1=13 D\else\ifnum#1=14 E\else\ifnum#1=15 F\fi\fi\fi\fi\fi\fi\fi}

\def\msa@{\hexnumber@\msafam}
\def\msb@{\hexnumber@\msbfam}
\mathchardef\boxdot="2\msa@00
\mathchardef\boxplus="2\msa@01
\mathchardef\boxtimes="2\msa@02
\mathchardef\square="0\msa@03
\mathchardef\blacksquare="0\msa@04
\mathchardef\centerdot="2\msa@05
\mathchardef\lozenge="0\msa@06
\mathchardef\blacklozenge="0\msa@07
\mathchardef\circlearrowright="3\msa@08
\mathchardef\circlearrowleft="3\msa@09
\mathchardef\rightleftharpoons="3\msa@0A
\mathchardef\leftrightharpoons="3\msa@0B
\mathchardef\boxminus="2\msa@0C
\mathchardef\Vdash="3\msa@0D
\mathchardef\Vvdash="3\msa@0E
\mathchardef\vDash="3\msa@0F
\mathchardef\twoheadrightarrow="3\msa@10
\mathchardef\twoheadleftarrow="3\msa@11
\mathchardef\leftleftarrows="3\msa@12
\mathchardef\rightrightarrows="3\msa@13
\mathchardef\upuparrows="3\msa@14
\mathchardef\downdownarrows="3\msa@15
\mathchardef\upharpoonright="3\msa@16

\mathchardef\downharpoonright="3\msa@17
\mathchardef\upharpoonleft="3\msa@18
\mathchardef\downharpoonleft="3\msa@19
\mathchardef\rightarrowtail="3\msa@1A
\mathchardef\leftarrowtail="3\msa@1B
\mathchardef\leftrightarrows="3\msa@1C
\mathchardef\rightleftarrows="3\msa@1D
\mathchardef\Lsh="3\msa@1E
\mathchardef\Rsh="3\msa@1F
\mathchardef\rightsquigarrow="3\msa@20
\mathchardef\leftrightsquigarrow="3\msa@21
\mathchardef\looparrowleft="3\msa@22
\mathchardef\looparrowright="3\msa@23
\mathchardef\circeq="3\msa@24
\mathchardef\succsim="3\msa@25
\mathchardef\gtrsim="3\msa@26
\mathchardef\gtrapprox="3\msa@27
\mathchardef\multimap="3\msa@28
\mathchardef\therefore="3\msa@29
\mathchardef\because="3\msa@2A
\mathchardef\doteqdot="3\msa@2B

\mathchardef\triangleq="3\msa@2C
\mathchardef\precsim="3\msa@2D
\mathchardef\lesssim="3\msa@2E
\mathchardef\lessapprox="3\msa@2F
\mathchardef\eqslantless="3\msa@30
\mathchardef\eqslantgtr="3\msa@31
\mathchardef\curlyeqprec="3\msa@32
\mathchardef\curlyeqsucc="3\msa@33
\mathchardef\preccurlyeq="3\msa@34
\mathchardef\leqq="3\msa@35
\mathchardef\leqslant="3\msa@36
\mathchardef\lessgtr="3\msa@37
\mathchardef\backprime="0\msa@38
\mathchardef\risingdotseq="3\msa@3A
\mathchardef\fallingdotseq="3\msa@3B
\mathchardef\succcurlyeq="3\msa@3C
\mathchardef\geqq="3\msa@3D
\mathchardef\geqslant="3\msa@3E
\mathchardef\gtrless="3\msa@3F
\mathchardef\sqsubset="3\msa@40
\mathchardef\sqsupset="3\msa@41
\mathchardef\trianglerighteq="3\msa@44
\mathchardef\trianglelefteq="3\msa@45
\mathchardef\bigstar="0\msa@46
\mathchardef\between="3\msa@47
\mathchardef\blacktriangledown="0\msa@48
\mathchardef\blacktriangleright="3\msa@49
\mathchardef\blacktriangleleft="3\msa@4A
\mathchardef\blacktriangle="0\msa@4E
\mathchardef\triangledown="0\msa@4F
\mathchardef\eqcirc="3\msa@50
\mathchardef\lesseqgtr="3\msa@51
\mathchardef\gtreqless="3\msa@52
\mathchardef\lesseqqgtr="3\msa@53
\mathchardef\gtreqqless="3\msa@54
\mathchardef\Rrightarrow="3\msa@56
\mathchardef\Lleftarrow="3\msa@57
\mathchardef\veebar="2\msa@59
\mathchardef\barwedge="2\msa@5A
\mathchardef\doublebarwedge="2\msa@5B
\mathchardef\angle="0\msa@5C
\mathchardef\measuredangle="0\msa@5D
\mathchardef\sphericalangle="0\msa@5E
\mathchardef\varpropto="3\msa@5F
\mathchardef\smallsmile="3\msa@60
\mathchardef\smallfrown="3\msa@61
\mathchardef\Subset="3\msa@62
\mathchardef\Supset="3\msa@63
\mathchardef\Cup="2\msa@64

\mathchardef\Cap="2\msa@65

\mathchardef\curlywedge="2\msa@66
\mathchardef\curlyvee="2\msa@67
\mathchardef\leftthreetimes="2\msa@68
\mathchardef\rightthreetimes="2\msa@69
\mathchardef\subseteqq="3\msa@6A
\mathchardef\supseteqq="3\msa@6B
\mathchardef\bumpeq="3\msa@6C
\mathchardef\Bumpeq="3\msa@6D
\mathchardef\lll="3\msa@6E

\mathchardef\ggg="3\msa@6F

\mathchardef\circledS="0\msa@73
\mathchardef\pitchfork="3\msa@74
\mathchardef\dotplus="2\msa@75
\mathchardef\backsim="3\msa@76
\mathchardef\backsimeq="3\msa@77
\mathchardef\complement="0\msa@7B
\mathchardef\intercal="2\msa@7C
\mathchardef\circledcirc="2\msa@7D
\mathchardef\circledast="2\msa@7E
\mathchardef\circleddash="2\msa@7F
\def\ulcorner{\delimiter"4\msa@70\msa@70 }
\def\urcorner{\delimiter"5\msa@71\msa@71 }
\def\llcorner{\delimiter"4\msa@78\msa@78 }
\def\lrcorner{\delimiter"5\msa@79\msa@79 }
\def\yen{\mathhexbox\msa@55 }
\def\checkmark{\mathhexbox\msa@58 }
\def\circledR{\mathhexbox\msa@72 }
\def\maltese{\mathhexbox\msa@7A }
\mathchardef\lvertneqq="3\msb@00
\mathchardef\gvertneqq="3\msb@01
\mathchardef\nleq="3\msb@02
\mathchardef\ngeq="3\msb@03
\mathchardef\nless="3\msb@04
\mathchardef\ngtr="3\msb@05
\mathchardef\nprec="3\msb@06
\mathchardef\nsucc="3\msb@07
\mathchardef\lneqq="3\msb@08
\mathchardef\gneqq="3\msb@09
\mathchardef\nleqslant="3\msb@0A
\mathchardef\ngeqslant="3\msb@0B
\mathchardef\lneq="3\msb@0C
\mathchardef\gneq="3\msb@0D
\mathchardef\npreceq="3\msb@0E
\mathchardef\nsucceq="3\msb@0F
\mathchardef\precnsim="3\msb@10
\mathchardef\succnsim="3\msb@11
\mathchardef\lnsim="3\msb@12
\mathchardef\gnsim="3\msb@13
\mathchardef\nleqq="3\msb@14
\mathchardef\ngeqq="3\msb@15
\mathchardef\precneqq="3\msb@16
\mathchardef\succneqq="3\msb@17
\mathchardef\precnapprox="3\msb@18
\mathchardef\succnapprox="3\msb@19
\mathchardef\lnapprox="3\msb@1A
\mathchardef\gnapprox="3\msb@1B
\mathchardef\nsim="3\msb@1C
\mathchardef\napprox="3\msb@1D
\mathchardef\nsubseteqq="3\msb@22
\mathchardef\nsupseteqq="3\msb@23
\mathchardef\subsetneqq="3\msb@24
\mathchardef\supsetneqq="3\msb@25
\mathchardef\subsetneq="3\msb@28
\mathchardef\supsetneq="3\msb@29
\mathchardef\nsubseteq="3\msb@2A
\mathchardef\nsupseteq="3\msb@2B
\mathchardef\nparallel="3\msb@2C
\mathchardef\nmid="3\msb@2D
\mathchardef\nshortmid="3\msb@2E
\mathchardef\nshortparallel="3\msb@2F
\mathchardef\nvdash="3\msb@30
\mathchardef\nVdash="3\msb@31
\mathchardef\nvDash="3\msb@32
\mathchardef\nVDash="3\msb@33
\mathchardef\ntrianglerighteq="3\msb@34
\mathchardef\ntrianglelefteq="3\msb@35
\mathchardef\ntriangleleft="3\msb@36
\mathchardef\ntriangleright="3\msb@37
\mathchardef\nleftarrow="3\msb@38
\mathchardef\nrightarrow="3\msb@39
\mathchardef\nLeftarrow="3\msb@3A
\mathchardef\nRightarrow="3\msb@3B
\mathchardef\nLeftrightarrow="3\msb@3C
\mathchardef\nleftrightarrow="3\msb@3D
\mathchardef\divideontimes="2\msb@3E
\mathchardef\varnothing="0\msb@3F
\mathchardef\nexists="0\msb@40
\mathchardef\mho="0\msb@66
\mathchardef\thorn="0\msb@67
\mathchardef\beth="0\msb@69
\mathchardef\gimel="0\msb@6A
\mathchardef\daleth="0\msb@6B
\mathchardef\lessdot="3\msb@6C
\mathchardef\gtrdot="3\msb@6D
\mathchardef\ltimes="2\msb@6E
\mathchardef\rtimes="2\msb@6F
\mathchardef\shortmid="3\msb@70
\mathchardef\shortparallel="3\msb@71
\mathchardef\smallsetminus="2\msb@72
\mathchardef\thicksim="3\msb@73
\mathchardef\thickapprox="3\msb@74
\mathchardef\approxeq="3\msb@75
\mathchardef\succapprox="3\msb@76
\mathchardef\precapprox="3\msb@77
\mathchardef\curvearrowleft="3\msb@78
\mathchardef\curvearrowright="3\msb@79
\mathchardef\digamma="0\msb@7A
\mathchardef\varkappa="0\msb@7B
\mathchardef\hslash="0\msb@7D
\mathchardef\hbar="0\msb@7E
\mathchardef\backepsilon="3\msb@7F
\def\Bbb{\ifmmode\let\next\Bbb@\else
 \def\next{\errmessage{Use \string\Bbb\space only in math mode}}\fi\next}
\def\Bbb@#1{{\Bbb@@{#1}}}
\def\Bbb@@#1{\fam\msbfam#1}

\catcode`\@=\active

\def\inv{^{\raise.15ex\hbox{${
  \scriptscriptstyle -}$}\kern-.05em 1}}

\def\Dsl{\,\raise.15ex\hbox{$/$}\mkern-13.5mu D}
\def\dsl{\raise.15ex\hbox{$/$}\kern-.57em\hbox{$\partial$}}

\def\lspace{\ifx\answ\bigans{}\else\qquad\fi}

\def\CL{\hbox{{$\cal L$}}} \def\CH{\hbox{{$\cal H$}}}

\def\CR{\hbox{{$\cal R$}}}

\def\lform{\hbox{$\sqcup$}\llap{\hbox{$\sqcap$}}}
\def\darr#1{\raise1.5ex\hbox{$\leftrightarrow$}
\mkern-16.5mu #1}

\def\INT{{\textstyle \int\kern-.642em\int}}

\def\R{{\Bbb R}}

\def\eps{{\epsilon}}

\def\rcross{{\triangleright\!\!\!<}}
\def\lcross{{>\!\!\!\triangleleft}}

\def\rbiprod{{\cdot\kern-.33em\triangleright\!\!\!<}}
\def\lbiprod{{>\!\!\!\triangleleft\kern-.33em\cdot}}

\def\tens{\mathop{\otimes}}

\def\la{{\triangleright}}\def\ra{{\triangleleft}}
\def\cora{{\blacktriangleright}}

\def\swap{{\leftrightarrow}}

\def\isom{{\cong}}

\def\span{{\rm span}}

\def\Ad{{\rm Ad}}

\def\id{{\rm id}}

\def\<{\langle}
\def\>{\rangle}
\def\dila{{\varsigma}}

\def\haj#1{{\mathaccent20 {#1}}}

\def\vect{{\bf t}}\def\vecs{{\bf s}}
\def\vecu{{\bf u}}

\def\<{\langle}
\def\>{\rangle}

\def\equad{\kern -1.7em}
\def\qqquad{\qquad\quad}

\def\o{{}_{\scriptscriptstyle(1)}}
\def\t{{}_{\scriptscriptstyle(2)}}
\def\th{{}_{\scriptscriptstyle(3)}}
\def\fo{{}_{\scriptscriptstyle(4)}}

\def\bo{{}^{\bar{\scriptscriptstyle(1)}}}
\def\bt{{}^{\bar{\scriptscriptstyle(2)}}}

\def\und#1{{\underline {#1}}}

\def\uo{{{}^{\scriptscriptstyle(1)}}}

\def\umo{{{}^{\scriptscriptstyle-(1)}}}
\def\umt{{{}^{\scriptscriptstyle-(2)}}}

\def\Bo{{{}_{\und{\scriptscriptstyle(1)}}}}
\def\Bt{{{}_{\und{\scriptscriptstyle(2)}}}}

\def\text#1{\mbox{\rm #1}}
\def\note#1{}

\def\blacksquare{{\lform}}
\def\frac#1#2{{{#1\over#2}}}

\def\proof{\goodbreak\noindent{\bf Proof\quad}}

\def\endproof{{\ $\lform$}\bigskip }

\def\eqn#1#2{\begin{equation}#2\label{#1}\end{equation}}

\def\align#1{\begin{eqnarray*}#1\end{eqnarray*}}

\def\ceqn#1#2{\begin{equation}\label{#1}\begin{array}{c}#2\end{array}
\end{equation}}

\documentstyle[11pt, epsf]{article}
\newtheorem{lemma}{Lemma}[section] \newtheorem{propos}[lemma]{Proposition}

\begin{document}\baselineskip 22pt

{\ }\qquad\qquad Submitted Proc. Banach Center Minisemester on QG, November 1995 
\vspace{.2in}

\begin{center} {\LARGE SOME REMARKS ON QUANTUM AND BRAIDED\\ GROUP GAUGE THEORY}
\\ \baselineskip 13pt{\ }
{\ }\\
 S. Majid\footnote{Royal Society University Research Fellow and Fellow of
Pembroke College, Cambridge}\\
{\ }\\
Department of Mathematics, Harvard University\\
Science Center, Cambridge MA 02138, USA\footnote{During 1995+1996}\\
+\\
Department of Applied Mathematics \& Theoretical Physics\\
University of Cambridge, Cambridge CB3 9EW\\
\end{center}
 
\vspace{10pt}
\begin{quote}
\baselineskip 13pt
\noindent{\bf Abstract}  We clarify some aspects of quantum group gauge theory and  its recent generalisations (by  T. Brzezinski and the author) to braided group gauge theory and coalgebra gauge theory. We outline the diagrammatic version of the braided case. We study the bosonisation of any braided group provides as a trivial principal bundle in three ways.
 
\bigskip

\end{quote}
\baselineskip 15pt

\section{Introduction} 

Quantum group gauge theory over quantum spaces has been introduced some years ago by T. Brzezinski and the author in \cite{BrzMa:gau}, which contained not only the formalism (which is fairly straightforward) but the non-trivial example of the  $q$-monopole to justify it. In spite of this success and a couple of follow-on papers \cite{BrzMa:mon}\cite{Brz:tra}\cite{Haj:str} extending the formalism of \cite{BrzMa:gau}, the $q$-monopole remains the main topologically non-trivial example of the formalism (over a quantum base). It even appears that for the $q$-analogue of other non-trivial bundles in physics  one may need to develop some generalisation beyond the quantum group case. 

In fact, just such  a generalised gauge theory has been recently introduced in \cite{BrzMa:coa}. The theory has for its `structure group' merely a coalgebra. The total space of the principal bundle and the base space are both algebras, possibly non-commutative. None of these structures need be quantum groups of any kind and hence it is remarkable that a full gauge theory (with bundles, connections, gauge transformations etc) is still possible. See Brzezinski's own contribution to these proceedings for an announcement. The theory  is general enough  to include the `embeddable homogeneous spaces' of \cite{Brz:hom}. Moreover, another special case of the theory is {\em braided group gauge theory} with structure group a braided group, introduced by these means in \cite{BrzMa:coa}. 

In this note we collect some modest results and remarks on the quantum group gauge theory of \cite{BrzMa:gau} and these generalisations. In the case of \cite{BrzMa:gau} we clarify in Section~2 a few points of formalism in the light of subsequent works. In Section~3 we announce a fully braid-diagrammatic version of the braided group gauge theory\cite{Ma:dia}. In Section~4 we point out that the inhomogeneous quantum groups $H\rcross B$ obtained by bosonisation of braided groups (such as $q$-Poincar\'e groups\cite{Ma:poi}) can be viewed in {\em all three} ways, as quantum principal bundles, as braided principal bundles, and as embeddable homogeneous spaces.

\subsection*{Acknowledgements}

I would like to thank T. Brzezinski for our continuing discussions under EPSRC research grant GR/K02244, from   which this work is an outgrowth. I also want to thank the organisers R. Budzynski, W. Pusz and S.~Zakrewski for a superb conference and minisemester at the Banach Center in Warsawa during November-December, 1995.

\section{Quantum group gauge theory}

Our comments here are confined to the approach initiated in \cite{BrzMa:gau}. In particular, we would not feel competent to comment on the interesting variant of this approach which was described at this conference by Mico Durdevic. Much of his formalism is the same as \cite{BrzMa:gau} but there are interesting differences too.  

We recall first the formalism of \cite{BrzMa:gau}. We take for the fiber of a principal bundle a quantum group $H$. The total space of the bundle is an $H$ comodule-algebra $P$, i.e. there is an algebra map $P\to P\tens H$, denoted $\cora(u)= u\bo\tens u\bt$, forming a coaction. We use notations as in the author's text \cite{Ma:book} and we denote by $M=\{u\in P|\cora(u)=u\tens 1\}$ the fixed point subalgebra of $P$, which is the base space of the bundle. We assume that $P$ is flat as an $M$-module. The composite $\tilde\chi:P\tens P\to P\tens B$ sending $u\tens v\to uv\bo\tens v\bt$ (denoted $\tilde{\ }$ in \cite{BrzMa:gau}) plays the role of the left-invariant vector fields generated by the group action.  We require it to be surjective, corresponding to freeness of the action. Finally, we require injectivity of the induced map $\chi:P\tens_M P\to P\tens B$, i.e. we require $\chi$ to be an isomorphism (the Galois condition). In \cite{BrzMa:gau} we imposed the injectivity of $\chi$ in the form $\ker\tilde\chi=P(\Omega^1M)P$ (a `differential Galois condition' ) which seemed more natural at the time but turned out to be equivalent. In fact, there still remain some problems with the correct gauge-invariant formulation of `smoothness' for the trivialisation  for non-universal differential calculus and one might want to come back to something more like the differential form of the injectivity at a future point. Concerning differential calculi,  we use the notation in \cite{BrzMa:gau}. We concentrate for simplicity on the universal calculus but one can fill in the nonuniversal case along the lines in \cite{BrzMa:gau}.

The maps $\chi$ and $\tilde\chi$ play a central role in \cite{BrzMa:gau}. From the axioms of a comodule algebra one sees at once their covariance properties. Namely, as intertwiners  
\eqn{covchi}{  \tilde\chi:P\tens P_\cora\to P\tens H_R,\quad  \tilde\chi:P_\cora\tens P\to P_\cora\tens H_L,\quad \Rightarrow \quad \tilde\chi:P_\cora\tens P_\cora\to P_\cora\tens H_\Ad}
and similarly for $\chi$. Here $P_\cora$ denotes $P$ taken with its given coaction $\cora$,  $H_R$ denotes $H$ with the regular right coaction $\Delta$, $H_L$ with the right-coaction version $h\mapsto h\t\tens Sh\o$ of the left regular representation and $H_\Ad$ the right adjoint coaction.  
The first two properties are elementary calculations used in the proof of the third Ad-invariance property.  This, in turn, is central to the construction in \cite{BrzMa:gau} of a connection $\Pi$ from a connection form $\omega$. 

The converse direction is also essentially in \cite{BrzMa:gau}, but no explicit formula for $\omega$ in terms of $\Pi$.  We provide this now as (as far as I know) a modest new result. We recall that a connection form is an intertwiner $H_\Ad\to \Omega^1P$ such that $\tilde\chi\circ\omega(u)=1-\eps(u)$ and $\omega(1)=0$. Note that $H=\{1\}\oplus\ker\eps$ so fixing the value of $\omega$ on 1 (as preferred in \cite{BrzMa:gau}) is equivalent to only specifying it on $\ker\eps$

\begin{propos} Let $P,\cora$ be a quantum principal bundle. There is a 1-1 correspondence between connection forms $\omega:H_\Ad\to \Omega^1P$ and covariant projections $\Pi:\Omega^1P\to\Omega^1P$  which are left $P$-module maps and have $\ker\Pi=P(\Omega^1M)P$. Explicitly,
\eqn{Piomega}{\Pi=(\cdot\tens\id)\circ(\id\tens\omega)\circ\tilde\chi}
as in \cite{BrzMa:gau}, and
\eqn{omegaPi}{ \omega(h)=\Pi\circ\chi^{-1}(1\tens h),\quad \forall h\in \ker\eps.}
\end{propos}
\proof Only the converse direction from $\Pi$ to $\omega$ needs to be clarified here. Indeed, we can obviously rewrite (\ref{covchi}) at the level of $\chi$ (without any calculation) as
\eqn{covchiinv}{ \chi^{-1}:P\tens H_R\to P\tens_M P_\cora,\quad   \chi^{-1}:P_\cora\tens H_L\to P_\cora\tens_M P,\quad \Rightarrow \quad  \chi^{-1}:P_\cora\tens H_\Ad\to P_\cora\tens_M P_\cora.}
Then given $\Pi$, we define $\omega$ via (\ref{omegaPi}). The formula makes sense because $P(dM)P\subset P(\Omega^1M)P\subseteq \ker\Pi$, where $dm=1\tens m-m\tens 1$ for $m\in M$, i.e. $\Pi$ descends to   $P\tens_M P$. We also need for the formula to make sense that $\cdot\circ\chi^{-1}(1\tens h)=(\id\tens\eps)\chi\circ\chi^{-1}(1\tens h)=\eps(h)=0$ when $h\in \ker\eps$, since $\Pi$ is defined on $\Omega^1P=\ker\cdot$. The intertwining condition for $\omega$ follows from (\ref{covchiinv}) and  $\tilde\chi\circ\omega(h)=\tilde\chi\circ\chi^{-1}(1\tens h)=1\tens h$, where   $\tilde\chi\circ\Pi=\tilde\chi$ follows from $\ker\tilde\chi=P(\Omega^1M)P\supseteq\ker\Pi$ and $\Pi^2=\Pi$. So we have a connection. Starting with a connection, we can consider the corresponding projection  
$\Pi_\omega$ from (\ref{Piomega}) and apply the above formula (\ref{omegaPi}) -- we clearly get back $\omega$. In the other direction, let $\omega$ be defined from $\Pi$. Then $\Pi_\omega(u\tens v)=uv\bo\Pi\circ\chi^{-1}(1\tens v\bt)=\Pi\left(u v\bo\chi\umo(1\tens v\bt)\tens \chi\umt(1\tens v\bt)\right)=\Pi\left(u\tens v\bo\chi\umo(1\tens v\bt) \chi\umt(1\tens v\bt)\right)=\Pi(u\tens v)$, where   $v\bo\chi\umo(1\tens v\bt)\in M$ so that we can move it over in $P\tens_MP$. This last fact requires two steps. The first is to compute the coaction on the first factor of $v\bo\chi\umo(1\tens v\bt)\tens_M\chi^{-1}(1\tens v\bt)$ and see that it is  trivial. This is immediate from the middle property in (\ref{covchiinv}) and is identical to the proof in \cite[p. 604]{BrzMa:gau} that $u\bo\Phi^{-1}(u\bt\o)\in M$ when $\Phi^{-1}:H_L\to P_\cora$. Indeed, $\chi^{-1}$ plays the role for a general bundle that $\Phi^{-1}$ plays for a trivial bundle. The only difference, which is the new feature of the present proof, is that from the trivial coaction on the first factor of an element $\sum u_i\tens_Mv_i$ we cannot in general conclude that $u_i\in M$. For this we need that $(\ )\tens_MP$ is a left exact functor, i.e. that $P$ is flat as an $M$-module. This is the case in most examples, however. \endproof

Another point which should be clarified from \cite{BrzMa:gau} is that the gauge transformations introduced there, which are convolution-invertible unit-preserving maps $\gamma:H\to M$, are only relevant for trivial bundles, i.e. they are really part of the local theory. This is the main setting in which one wants gauge transformations (to patch trivial bundles) and was emphasised for this reason. On the other hand, just as gauge fields $A:H\to \Omega^1M$ on trivial bundles lead to globally-defined connection forms $\omega:H_\Ad\to\Omega^1P_\cora$ by the formula
\eqn{omegaA}{  \omega_{A,P,\Phi}=\Phi^{-1}*d\Phi+\Phi^{-1}*A*\Phi}
(where $*$ denotes the convolution product and $\Phi$ the trivialisation) it is obvious that local gauge transformations $\gamma:H\to M$ lead to global gauge transformations $\Gamma:H_\Ad\to P_\cora$ by
\eqn{Gamgam}{ \Gamma_{\gamma,P,\Phi}=\Phi^{-1}*\gamma*\Phi}
Indeed, the calculations are identical to (a simpler version of) those for $\omega$ in \cite{BrzMa:gau} and we do not need to repeat them. The global $\omega,\Gamma$ can then be formulated for general bundles. This is the strategy developed in \cite{BrzMa:gau} and mirrors the strategy  in classical geometry whereby pseudotensorial forms on $P$ are the correct generalisation in classical geometry of local fields on the base. The correspondence for general  pseudotensorial forms is also treated in \cite{BrzMa:gau} in the quantum case, in addition to the above adjoint representation which has to be treated specially via $\Phi^{-1}*(\ )*\Phi$, as explained in \cite{BrzMa:gau}. A variant of \cite[Prop.~A.7]{BrzMa:gau} checks that $\gamma=\Phi*\Gamma*\Phi^{-1}$ has its values in $M$.

In an addendum\cite{Brz:tra} to  \cite{BrzMa:gau}, T. Brzezinski has taken the natural step to define from a general gauge transform (a convolution-invertible unit-preserving intertwiner) $\Gamma:H_\Ad\to P_\cora$ a so-called `bundle automorphism' $\Theta:P_\cora\to P_\cora$ defined by \eqn{ThetaGam}{\Theta(u)=u\bo\Gamma(u\bt).} 
The formula is similar to the way that $\omega$ defines a covariant  $\Pi$ and by a similar computation  one has that $\Theta$ is covariant (an intertwiner for the coaction). It is also immediate from the comodule algebra axioms that $\Theta$ is a left $M$-module map and that $\Theta(1)=1$ and $\Theta_\Gamma\circ\Theta_{\Gamma'}=\Theta_{\Gamma'*\Gamma}$ (so that $\Theta$ is invertible). While these steps are all immediate, Brzezinski succeeded to explicitly prove the converse (which is harder), namely that from $\Theta:P\to P$ obeying these properties one can define\cite{Brz:tra}
\eqn{GamTheta}{ \Gamma(h)=\chi\umo(1\tens h)\Theta(\chi\umt(1\tens h)).}
Brzezinski's proof is one of the motivations behind our
proof of Proposition~2.1 above. The main differences in the proof are that we require flatness of $P$ and we deduce the covariance properties of $\chi^{-1}$ from those of $\chi$ without proving them directly as is done in \cite{Brz:tra} under the heading of properties of the `translation map'  $\chi^{-1}(1\tens( ))$. Modulo these differences, we see that Proposition~2.1 completes the analogy between $\Pi,\omega,A$ and $\Theta,\Gamma,\gamma$ by providing the correspondence between $\Pi$ and $\omega$ in the same form as between $\Theta$ and $\Gamma$. 

For our second modest result we clarify further that $\Theta$ is not really a  bundle automorphisms in a natural way, because it is  not an algebra map. Instead, we propose the following definition which is natural from the point of view of extension theory: Fix  $M$ and consider all principal bundles with base $M$. Two of them are equivalent if there is a comodule algebra isomorphism between them which preserves $M$. In this sense, we should really view $\Theta$ as a {\em bundle gauge transformations} $\Theta:P\to P^\Gamma$, where $P^\Gamma=P$ as an $H$-comodule but has a new product 
\eqn{PGam}{ u\cdot_\Gamma v=\Theta(\Theta^{-1}(u)\Theta^{-1}(v))=u\bo\Gamma^{-1}(u\bt\o)v\bo\Gamma^{-1}(v\bo\o)\Gamma(u\bt\t v\bt\t).}
This also forms a bundle with the same base $M$ and $\Theta:P\to P^\Gamma$ is a bundle equivalence transforming our original bundle to it. This is an interesting {\em new phenomenon} in non-commutative geometry because we see from (\ref{PGam}) that $\cdot_\Gamma=\cdot$ when our algebras are commutative and $\Gamma$ are restricted to algebra maps. This is not possible in the general non-commutative case and so we have this new phenomenon that {\em the bundle itself gauge transforms!}

\begin{propos} Let $\Theta:P\to P^\Gamma$ be a bundle gauge transformation. If $\omega$ is a connection on $P$ with projection $\Pi$ then $\omega^\Gamma=(\Theta\tens\Theta)\circ\omega$ is a connection on $P^\Gamma$ and its associated projection $\Pi^\Gamma$ obeys $(\Theta\tens\Theta)\circ\Pi=\Pi^\Gamma\circ(\Theta\tens\Theta)$. Moreover, if $\Phi$ is a trivialisation of $P$ then $\Phi^\Gamma=\Theta\circ\Phi= \Phi*\Gamma$ is a trivialisation of $P^\Gamma$ such that  
\[ (\omega_{A,P,\Phi})^\Gamma=\omega_{A,P^\Gamma,\Phi^\Gamma}\]
\end{propos}
\proof This is elementary; we apply $\Theta$ systematically to all constructions in $P$. Then by definition all constructions in $P^\Gamma$ are the same after allowing for this algebra isomorphism $\Theta:P
\isom P^\Gamma$. Because $\Theta$ is an $H$-comodule map and $M$-module map, we do not need to change the coaction $\cora$ or $M$. For example, $\Pi^\Gamma(u\tens v)=u\cdot_\Gamma v\bo\cdot_\Gamma \omega^\Gamma(v\bt)=(\Theta\tens\Theta)\left(\Theta^{-1}(u)\Theta^{-1}(v)\bo\omega(\Theta^{-1}(v)\bt)\right)=(\Theta\tens\Theta)(\Pi(\Theta(u)\tens\Theta(v))$
for $u,v\in P$. Likewise, when we compute (\ref{omegaA}) using the $\Phi^\Gamma$ and the product in $P^\Gamma$ we obtain $\omega^\Gamma$.  
Finally, $\Phi^\Gamma(h)=\Phi(h)\bo\Gamma(\Phi(h)\bt)=\Phi(h\o)\Gamma(h\t)$ as stated, due to the covariance of $\Phi$. \endproof

When $\Gamma$ is obtained from $\gamma$ via (\ref{Gamgam}) we have $\Phi^\Gamma=\Phi*\Phi^{-1}*\gamma*\Phi=\gamma*\Phi=\Phi^\gamma$ and 
\eqn{Agam}{(\omega_{A,P,\Phi})^\Gamma=\omega_{A,P^\Gamma,\Phi^\gamma}=\omega_{A^\gamma,P^\Gamma,\Phi};\quad A^\gamma=\gamma^{-1}*A*\gamma+\gamma^{-1}*d\gamma.}
Apart from the additional transformation from $P$ to $P^\Gamma$ (which is not present in the classical situation) we see that the global gauge transformations induce the local transformation picture in \cite{BrzMa:gau}.
This ties up the global and local pictures in a way that we have not seen elsewhere.
 
Moreover, it was shown in \cite{BrzMa:gau} that every trivial principal bundle   has a canonical identification with the vector space $M\tens H$. As discussed more explicitly in \cite{Ma:clau}, its product then becomes that of  a cocycle cross product $M{}_c\lcross H$ for some convolution-invertible cocycle   $c:H\tens H\to M$ and cocycle-action $\alpha:H\tens M\to M$. This is the canonical form for a trivial principal bundle. Cocycle cross products are a standard algebraic construction and we refer to \cite[Chapter~6.3]{Ma:book} for the required formulae and conventions.

\begin{propos} If $M_c\lcross H$ is a cocycle cross product trivial bundle then 
$(M_c\lcross H)^{\gamma^{-1}}=M_{c^\gamma}\lcross H$ where $c^\gamma,\alpha^\gamma$ are the cocycle data  transformed by $\gamma$.
\end{propos}
\proof We first compute $\Theta$ associated to $\gamma$ in the case of $P=M_c\lcross H$. Here $\cora(m\tens h)=m\tens h\o\tens h\t$ is the coaction and $\Phi(h)=1\tens h$ and $\Phi^{-1}(h)=c^{-1}(Sh\t\tens h\th)\tens Sh\o$ is the trivialisation (see \cite[Prop.~6.3.6]{Ma:book}). In fact, we can avoid the explicit form of $\Phi^{-1}$ by using associativity of the product of $P$ and the special form of the coaction. Thus,
\align{\Theta(m\tens h)\equad &&=(m\tens h\o)\Gamma(h\t)=(m\tens h\o) \Phi^{-1}(h\t)(\gamma(h\th)\tens h\fo)\\
&&=(m\tens 1)\Phi(h\o)\Phi^{-1}(h\t)(\gamma(h\th)\tens h\fo)=m\gamma(h\o)\tens h\t}
This is of just the form  which  implements  $M_{c^\gamma}\lcross H\isom M_c\lcross\gamma$ in \cite[Prop. 6.3.5]{Ma:book}. \endproof

This ties up our quantum geometrical notion of  gauge transformation of a bundle with the algebraic theory of equivalence of cocycle cross products\cite{Doi:equ}; see \cite[Prop.~6.3.4]{Ma:book} for the required transformations $c^\gamma,\alpha^\gamma$ and a discussion of the corresponding non-Abelian cohomology $H^2(H,M)$. We see that this non-Abelian cohomology precisely classifies trivial bundles based on $M$ with structure $H$ up to bundle gauge transformation. Note that in the classical picture all cocycles are trivial due to commutativity of the algebras, so we do not see this phenomenon.

Next, we return to the construction of connections $A$ from gauge fields $\omega$. For general matter fields the global picture is to work on the total space, i.e. as intertwiners $V\to \Omega^nP_\cora$ (the  pseudotensorial forms). It was shown in \cite{BrzMa:gau} that not all of these come from the base on a trivial bundle  and the ones that do, and which are therefore the correct  way of working with matter fields on any bundle, are the {\em strongly tensorial} forms $V\to(\Omega^nM)P$. The subspace $(\Omega^nM)P$ was introduced for this purpose in \cite{BrzMa:gau}. Likewise, it was noted in \cite{BrzMa:gau} that not every connection $\omega:H_\Ad\to \Omega^1P_\cora$ is of the form (\ref{omegaA}) coming from a gauge field $A$ but an abstract characterisation of which $\omega$ are of this `strong' form was not given. This has recently been provided by P. Hajac in \cite{Haj:str}, where it is shown that these are precisely the connections obeying 
\eqn{somega}{ \cora(u\bo\omega(u\bt))=u\tens 1\tens 1-u\bo\tens 1\tens u\bt+u\bo\omega(u\bt)\tens 1.}
where $\cora$ acts on the product of $u\bo$ with the left factor in the output of $\omega$, and the $H$ part of the output of $\cora$ is placed to the far right. This is a kind of invariance condition for the left factor of the output of $\omega$ in place of requiring its values in $(\Omega^1M)P$. In fact, the condition (\ref{somega}) can be written more elegantly in terms of $\Pi$ as the condition $(\id-\Pi)du\in(\Omega^1M)P$ for all $u\in P$, see \cite{Haj:str}.

Another point of view is the following. It was explained in \cite{BrzMa:gau} that the connections $\omega_{A,P,\Phi}$ which come from the base have the property that the covariant derivative $D\Sigma=(\id-\Pi)d\Sigma$ on any strongly tensorial form $\Sigma:V\to (\Omega^nM)P_\cora$ is again strongly tensorial. This condition also makes sense globally and from a geometrical point of view one should define a connection as strong {\em iff} its covariant derivative sends all strongly horizontal forms to strongly horizontal forms. In one direction Hajac already observed that if $\omega$ is strong in the sense of (\ref{somega}) then $D\Sigma$ is strongly tensorial for any strongly tensorial $\Sigma$. The proof is analagous to that in \cite{BrzMa:gau} for trivial bundles; the computation for $D\Sigma$ reduces when $\Sigma:V\to (\Omega^nM)P$ to $(\id-\Pi)d$ on the rightmost output of $\Sigma$, i.e. to the 0-form case. The reduced 0-form case is then clear from Hajac's condition (\ref{somega}) in terms of $\Pi$. It is perhaps worth remarking that the converse is also true: If $D$ has this property  then take $\Sigma=\id:P_\cora\to P_\cora$, viewed as a strongly horizontal 0-form. 
Since  $D(\id)=(\id-\Pi)d\circ\id=(\id-\Pi)d$ we see that $D(\id)$ being strongly tensorial is exactly the $\Pi$ form of Hajac's condition (\ref{somega}). So these two definitions of a connection  strong coincide. All nontrivial results here are surely in \cite{Haj:str} but we did not find it formulated this way and would like to emphasise it explicitly: a connection is strong {\em iff} its covariant derivative sends strongly tensorial forms to strongly tensorial ones.   

Finally, we would like to discuss the notion of associated fiber bundle to a principal bundle and problems associated with it. Our main observation is that in \cite[Appendix~A]{BrzMa:gau} we assumed that $V$ is an $H^{\rm op}$-comodule algebra in order that we could define $E=(P\tens V)^H$ as an algebra. The algebra contains $M$ and we call $E$ the associated fiber bundle (or, in relevant examples,  vector bundle), over $M$ with fiber $V$.  A cross section of $E$ should be a map $\vecs:E\to M$ which is a left $M$-module map and is the identity on $M$ (or equivalently, obeys $\vecs(1)=1$). The left $M$-module structure is obviously needed (since classically we can multiply sections by functions on the base) but was accidentally omitted in \cite{BrzMa:gau}; it is needed for a correct 1-1 correspondence in \cite[Cor~A.8]{BrzMa:gau} between cross sections $\vecs:E\to M$ and unit-preserving maps $\sigma:V\to M$ in the case of a trivial associated vector bundle. Since such `local sections' (or `matter fields') $\sigma$  also correspond as mentioned above to pseudotensorial 0-forms $V\to P$ we see that  the latter correspond to cross section, which makes sense globally. This global correspondence is easily   verified\cite{Brz:tra} by replacing the role of $\Phi^{-1}$ for trivial bundles by $\chi^{-1}(1\tens (\ ))$ in  for a general bundle (as with the $\Pi,\omega$ and $\Gamma,\gamma$ correspondences above.) 

While one can certainly make variations of these constructions to try to improve the appearance of some of the formulae, we want to propose here something more radical: when one looks over the proofs one finds that one does not actually need $V$ or $E$ to be algebras. We can assume only that $V$ is a right $H$-comodule equipped with an element $1\in V$ which is fixed under the coaction. We still have $P\tens V$ as the tensor product comodule and define $E=(P\tens V)^H$. It remains that $M\tens 1\subset E$ and $E$ is a left $M$-module. We can still define cross sections as unit-preserving left $M$-module maps $\vecs:E\to M$. Even in this setup, the correspondence with pseudotensorial forms $V\to P$ holds. Explicitly, it is
\eqn{phis}{ \vecs(u\tens v)=u\Sigma(v),\quad \Sigma(v)=\chi\umo(1\tens S^{-1}v\bt)\vecs\left(\chi\umt(1\tens S^{-1}v\bt)\tens v\bo\right).}
We also have the result in \cite{BrzMa:gau} that when $P$ is trivial then any associated $E$ also has a trivialisation $\Phi_E:V\to E$ giving an isomorphism
$E\isom M\tens V$. Explicitly, $\Phi_E(v)=\Phi(S^{-1}v\bt)\tens v\bo$, in fact by the same formula as in \cite{BrzMa:gau} even though our input assumptions are slightly different. Note that we do not make any effort here to eliminate $S^{-1}$ since all quasitriangular and dual quasitriangular Hopf algebras have invertible antipodes.  

It is not known if this more straightforward formulation (where we do not use $H^{\rm op}$) has any real advantages over the original formulation in \cite{BrzMa:gau} or its variants. The point is that the most natural associated bundle in gauge theory (after $V=H_R$) should be $V=H_\Ad$. But this is {\em not} a comodule algebra so does not fit into the setting of \cite{BrzMa:gau}. We see that this is not a problem and we can proceed in any case, with $1\in H_\Ad$ as our distinguished fixed element. So, for example, maps $\Gamma:H_\Ad\to P_\cora$ as above indeed correspond to cross sections of the adjoint bundle $E=(P\tens H_\Ad)^H$, with an additional condition corresponding to convolution-invertibility.  

This point of view also leads naturally to search for a braided generalisation of quantum group gauge theory. For the braided adjoint coaction of all braided groups obtained by transmutation, such as $BG_q$ associated to compact simple Lie groups\cite{Ma:exa}, {\em is} a braided comodule algebra structure. I.e. one can conjugate braided matrices in a way that one cannot conjugate quantum matrices. This has been explained in \cite{Ma:lin} as one of the key features if braided groups. We outline this in the next section. 

We can also use $BG_q$ as our fiber even when $H=G_q$ the corresponding matrix quantum group, and use the above formalism. So $E=(P\tens BG_q)^{G_q}$ is a natural definition of the adjoint bundle. We can also take the braided-Lie algebra $\CL=\span\{u^i{}_j\}$ which is stable under the adjoint coaction $\cora(\vecu)=\vect^{-1}\vecu\vect$. So $E=(P\tens \CL)^{G_q}$ has a finite-dimensional braided-Lie algebra fiber. We can also take $E=(P\tens U_q(g))^{G_q}$,  or more generally $E=(P\tens H^*)^H$ with a suitable dual and right coadjoint coaction. It may be possible in this way to transform (or perhaps, transmute) our  above theory where gauge transformations act by convolution-inversion as in (\ref{Agam}) to one in which they act by something more like the quantum adjoint coaction.

\note{The covariance (\ref{covchiinv}) tells us that $\Gamma:H_\Ad\to P_\cora$ as required. Moreover, starting with $\Theta$ and defining $\Gamma$ from it, we have $\Theta_\Gamma(v)=v\bo\chi\umo(1\tens v\bt)\Theta(\chi\umt(1\tens v\bt))=\Theta\left(v\uo\chi\umo(1\tens v\bt)\chi\umt(1\tens v\bt)\right)=\Theta(v)$ by the same argument as in the proof of Proposition~2.1 that $v\uo\chi\umo(1\tens v\bt)\in M$.}

\section{Braided group gauge theory}

Here we make some remarks on the braided group version of the gauge theory. This has been initiated in \cite{BrzMa:coa} with the key definitions for bundles, connections and trivialisations as an example of a still more general coalgebra gauge theory. However, by using braid-diagrammatic methods for braided groups\cite{Ma:introp} it is possible to extend this to include {\em all} of the formalism of quantum group gauge theory as it was explained in  the last section. This is done in \cite{Ma:dia}. As usual in braided mathematics, one has to choose carefully between under and over crossings and be sure that things do not get tangled up. It should also be possible to push all results up to the coalgebra gauge theory level as well, which is work in progress.

For braided groups we use the formalism in \cite{Ma:introp}. See also Chapters~9.4 and~10 of \cite{Ma:book}. Suffice it to recall that for a concrete braided group $B$ the coproduct is a homomorphism $\Delta:B\to B\und\tens B$, where the latter is the braided tensor product algebra $(b\tens c)(a\tens d)=b\Psi(c\tens a)d$. We also require a counit $\eps$, antipode $S:B\to B$ and identity element $1$. More generally, we work in a braided tensor category with all maps written as morphisms and the identity as a morphism $\eta:\und 1\to B$, where $\und 1$ is the trivial object. We will need also suitable direct sums, kernels etc. so we keep in mind the category of representations over some background quantum group. On  the other hand, all constructions will be universal or diagrammatic, which means that they can be formulated quite generally. We recall that all morphisms are written pointing generally downwards and we write $\Psi=\epsfbox{braid.eps}, \Psi^{-1}=\epsfbox{braidinv.eps}$ and $\cdot=\epsfbox{prodfrag.eps}$, $\Delta=\epsfbox{deltafrag.eps}$. Other morphisms are written as nodes, and the unit object $\und 1$ is denoted by omission. Axioms, existence, effective techniques for working and applications are due to the author in \cite{Ma:bra}\cite{Ma:exa}\cite{Ma:bos}\cite{Ma:lin}\cite{Ma:poi}\cite{Ma:lie} etc; see the \cite{Ma:introp}. 
 
We define a braided principal bundle to be a right braided comodule algebra\cite{Ma:introp} such that the map $\chi$ is an isomorphism, i.e. an algebra $P$ and a comodule structure $\cora:P\to P\und\tens B$ which is an algebra homomorphism to the braided tensor product algebra, and invertibility of $\chi$ induced from this data. Here $M=P^B$ is the kernel of $\cora-\id\tens \eta:P\to P\tens B$. Using braid diagrams one easily sees that $M$ is an algebra and that $\tilde\chi$ descends to $P\tens_MP$ (which is defined in the usual way). We assume that $\tens_M P$ is left exact.  We also define the universal differential calculus in the same way as before; the braiding does not enter there. Finally, using the braiding again, we extend $\cora$ to tensor products of $P$ by the braided tensor product coaction\cite{Ma:introp}. 
 
We can then define the structures which we had before, now as intertwiners for the coaction of $B$ (as well as morphisms in the category, i.e. intertwiners for the background quantum group in the concrete case). Thus
\eqn{intmors}{ \omega:B_\Ad\to \Omega^1P_\cora,\quad \Pi:\Omega^1P_\cora\to \Omega^1P_\cora,\quad \Gamma:B_\Ad\to P_\cora,\quad \Theta:P\cora\to P_\cora}
with the various conditions in the definitions being the same as before with the obvious modifications to express everything in terms of morphism (such as $\tilde\chi\circ\omega=\eta\tens (\id-\eta\circ\eps)$ and $\omega\circ\eta=0$ in the definition for a connection form). $\Ad$ in (\ref{intmors}) denotes the braided adjoint coaction introduced in \cite{Ma:lin}. We have 1-1 correspondences $\Pi\swap\omega$ and $\Theta\swap\Gamma$ as in the preceding section.

Let $V$ be any $B$-comodule. We define pseudotensorial forms, strongly tensorial forms  and the covariant derivative $D=(\id-\Pi)\circ d$ on them in the same way as before. The only difference is that $\Sigma:V\to\Omega^nP_\cora$ is an intertwiner to the braided tensor product coaction on $\Omega^nP$ and we use a braided version of the strongness condition, namely
\eqn{brasconn}{ \epsfbox{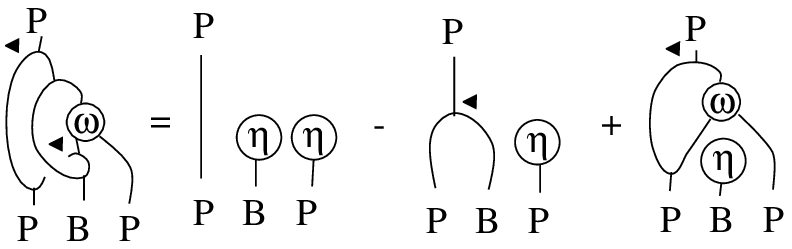}}
Then $D$ descends to strongly tensorial forms whenever the connection is strong.

Finally, suppose that $V$ is equipped with a morphism $\eta:\id\to V$ such that $\cora\circ\eta=\eta\tens\eta$. We say that $E=(P \tens V)^B$ is an associated braided fiber bundle with fiber $V$, where $P\und\tens V$ has the braided tensor product coaction. We have a left $M$-module structure on $E$ as before, and define a cross section as a $M$-module morphism $\vecs:E\to M$ such that $\vecs\circ\eta=\eta$. We have the correspondence $\vecs\swap\Sigma$ as before.

This describes the global theory. Next, we define a trivialisation as before, an intertwiner morphism $\Phi:B_R\to P_\cora$ which is convolution invertible and obeys $\Phi\circ\eta=\eta$. From such $\Phi$ on a braided comodule algebra $P$ one can make $P$ intro a trivial bundle. The isomorphisms $\theta:M\tens B\isom P$ and $\theta_E:M\tens V\isom E$ as objects in the category go through as in the quantum group case. 
 
Finally, we define local gauge fields, local gauge transformations and local sections or matter fields as morphisms
\eqn{localmors}{ A:B\to \Omega^1M,\quad \gamma:B\to M,\quad \sigma:V\to \Omega^nM}
as before. When $P$ is trivial the 1-1 correspondences strong $\omega\swap A$, $\Gamma\swap \gamma$ and $\Sigma\swap\sigma$ go through by just the same formulae as in the preceding section, expressed diagrammatically. 
\begin{figure}\[ \epsfbox{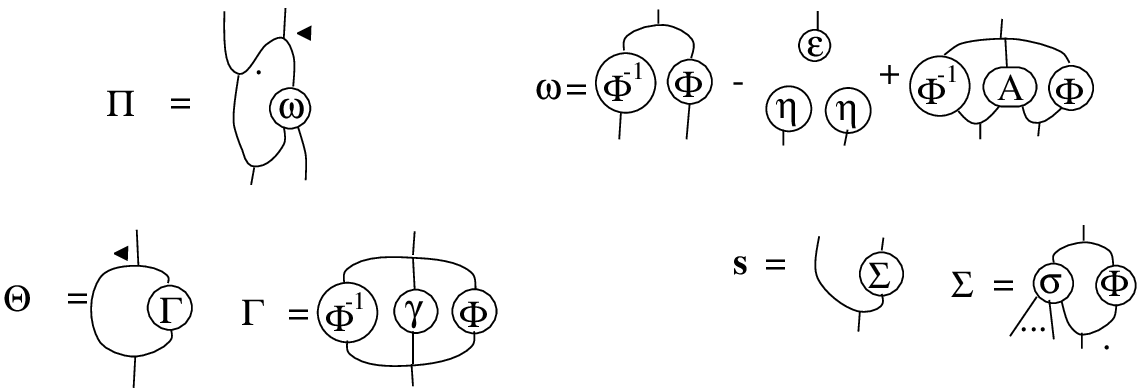}\]
\caption{Reduction from geometric level to local fields}
\end{figure}
 
The principal parts of the above correspondences are summarised in Figure~1. For the converses and other structures (and proofs) see \cite{Ma:dia}. The resulting picture at the level of these local fields works in just the same way as in the quantum group case in \cite{BrzMa:gau}. Thus, the covariant derivatives becomes $\nabla\sigma=d\sigma -(-1)^n\sigma*A$ and we have
\eqn{localF}{ \nabla^2\sigma=-\sigma*F, \quad F=dA+A*A.}
and the Bianchi identity $dF+A*F-F*A=0$ holds. The local gauge transformation induces gauge transformations $A^\gamma$ and $\sigma^\gamma$ and $\nabla^\gamma$, $F^\gamma$ computed with $A^\gamma$ are locally gauge covariant in the expected way. The diagrams behind the local theory are in Figure~2. Some of the $\epsfbox{deltafrag.eps}$ nodes refer to the coproduct and some to the coaction on $V$, as clear from context. We remark that {\em the local gauge theory does not require a braiding}. It works for a monoidal category with suitable direct sums. It does not use the product of $B$ either, i.e. needs only a coalgebra in the monoidal category.
\begin{figure}\[ \epsfbox{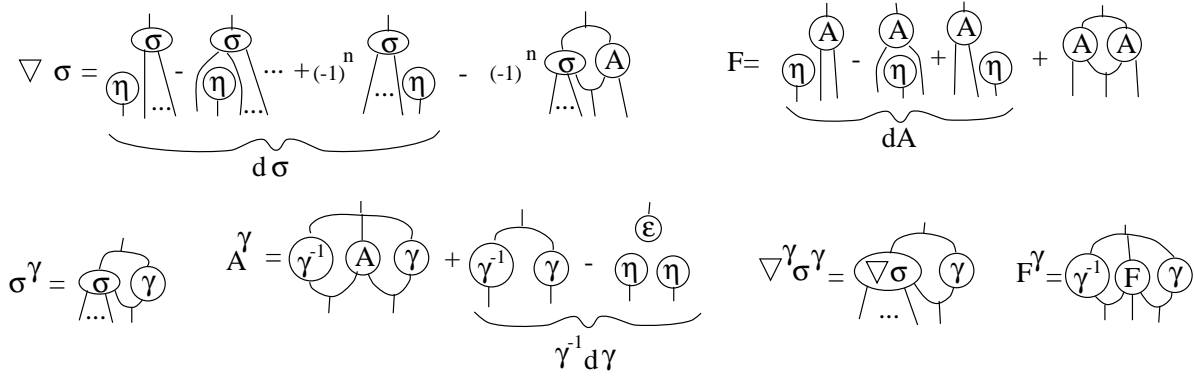}\]
\caption{Local form of braided gauge theory}
\end{figure}

We have obviously glossed over a lot here, but the basic message is that the entire theory generalises to the braided case provided the natural braided tensor product coactions and other structures are used. What is new in the braided case? Here are two key points. (i) The most trivial bundle is 
\eqn{tensbund}{ P=M\und\tens B,\quad \cora=\id\tens\Delta,\quad \Phi=\eta\tens\id,\quad 
\Phi^{-1}=\eta\tens S,}
but now we use the braided tensor product algebra, which can be complicated (the two factors do not commute). So we get complicated models very easily.
This is amply demonstrated in \cite{BrzMa:coa}. Moreover, $E=M\und\tens V$ in this case, again the braided tensor product algebra (when $V$ is a braided comodule algebra). So we obtain interesting associated bundles too. (ii) For braided groups obtained by full transmutation (such as $BG_q$) there is a natural braided commutativity condition
\eqn{bracom}{ \epsfbox{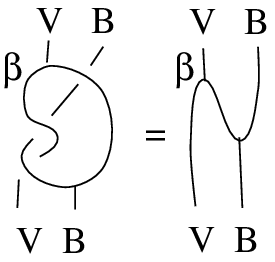}}
which holds  with respect to comodules $V$ (with coactions $\beta$) obtained by transmutation. When $B$ and a braided comodule algebra $V$ are of 
this type then $E=(P\und\tens V)^B$ becomes an algebra. Thus we restore the geometrical picture of $E$ as the `coordinate ring' on the total space without any unnatural use of $B^{\rm op}$ etc. In particular, for such braided groups we always have a natural adjoint bundle $E=(P\und\tens B_\Ad)^B$.

The completely dual version of the braided theory is easy enough: just turn all diagrams here and in \cite{Ma:dia} up side down! In this case our total space and base space are handled as coalgebras, which should be thought of as linear combinations of points, and $B$ becomes viewed as an enveloping algebra. A bit more involved is the semidualisation where we dualise only $B$. Using the dualisation  lemmas in \cite{Ma:introp} one obtains a theory where $P$ is a right $(B^{\haj{\ }})^{\rm cop}$ braided comodule algebra in the category with reversed braiding. Reflecting in a mirror and relabelling, we obtain the following formulation. $P$ a left braided $B$-module algebra with fixed subalgebra $M$, $\tilde\chi:B\tens P\tens P\to P$ (the box in Figure~3) descends to $\chi:B\tens P\tens_MP\to P$ and is non-degenerate in a certain sense. We denote the action of $P$ on $B$ by $\la$ and extend it to $\Omega^1P$. Connection forms, projections, 
global gauge transformations and bundle transformations are elements (more precisely, morphisms from $\und 1$) 
\eqn{envcon}{ \omega\in \Omega^1P_\la\tens B_\Ad,\quad\Pi:\Omega^1P_\la\to \Omega^1P_\la,\quad \Gamma\in P_\la\tens B_\Ad,\quad \Theta:P_\la\to P_\la}
which are also invariant under the action of $B$. Here $\Ad$ is the left braided adjoint action from \cite{Ma:lie} and we use the braided tensor product action on tensor product objects. The additional condition for $\omega$ is $\tilde\chi\circ\omega=\eta\tens\id -\eta\tens\eta\circ\eps$ and $\eps\circ\omega=0$, where the compositions are made in the way that make sense. The other parts of the definitions are likewise similar. Then we have 1-1 correspondences as before. Let $V$ be a right $B$-module by $\ra$. It is also useful to let $\la$ be the corresponding left module structure on $V$ defined via the braiding and braided-antipode as explained in\cite{Ma:introp}. Then define pseudotensorial forms as invariant elements $\Sigma\in V\tens P_\la$. We define associated vector bundles and sections much as before. A trivialisation is a suitably invertible such element $\Phi\in P_\la\tens B_L$ where $B_L$ is the left regular representation (it is actually the dualisation of $\Phi^{-1}$ in the comodule theory). On trivial bundles we obtain connections and gauge transformations from
\begin{figure}\[ \epsfbox{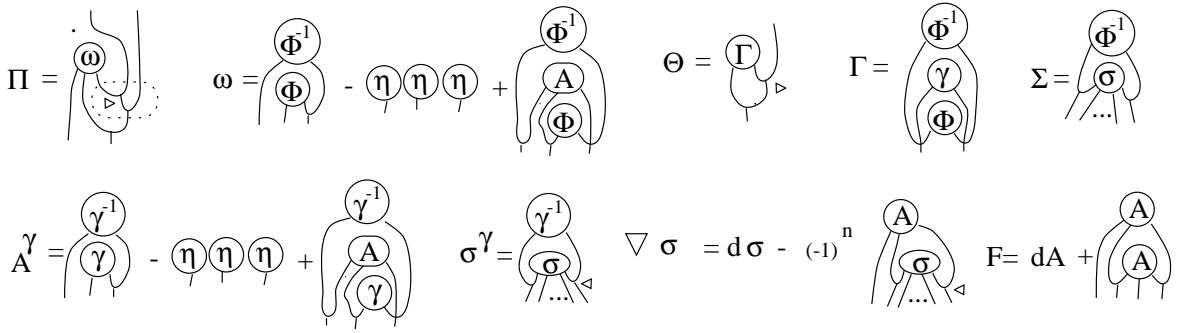}\]
\caption{Braided gauge theory of `enveloping algebra' type}
\end{figure}\eqn{envloccon}{ A\in  \Omega^1M\tens B,\quad \gamma\in M\tens B,\quad \sigma\in \Omega^nM\tens V.}
The key elements of the theory are in Figure~3. Note that convolution-invertibility in $M\tens B$ etc. becomes inverse in the `piggy back' product avoiding a braiding. This is the commuting tensor product $M\tens B^{\rm op}$ if one uses the usual opposite product. One can also reformulate everything in terms of the braided tensor product $M\und\tens B$ by a different dualisation involving the braiding (the natural categorical dual does not). In this case it is more natural to use the left action $\la$ in the diagrams involving $V$ in the local picture.

\section{Bosonisations as quantum/braided/homogeneous bundles}

There are plenty of braided algebras and groups, and hence plenty of braided   principal bundles. Here we want to focus on a a class of
braided tensor product  principal bundles which are closely connected with
quantum groups. They include the braided line example in \cite{BrzMa:coa} but also the extended $q$-Poincar\'e algebra studied extensively in \cite{Ma:qsta}.

Let $H$ be a dual-quasitriangular Hopf algebra with dual quasitriangular structure $\CR:H\tens H\to k$. We work in the braided category of $H$-comodules. Let $B$ be a braided group in this category. Bosonisation theory tells us that there is an ordinary Hopf algebra $H\rbiprod B$ given by the semidirect coproduct by the coaction of $H$ on $B$ as an object, and an action induced by $\CR$ (\cite{Ma:bos} in the dual form). Also from this theory, it is clear that
the bosonisation has the structure of a certain braided tensor product. Hence

\begin{propos} Every bosonisation in the braided category of right $H$-comodules can be viewed as  a   braided tensor product principal bundle 
$H\rbiprod B=H_R\und\tens B$ with structure braided group $B$ and the right coregular coaction $\cora=\id\tens\Delta$.
\end{propos}
\proof Let $H_R$ denote $H$ regarded as an $H$-comodule via the coproduct (the right coregular representation). Then the algebra structure of the bosonisation $H\rbiprod B$ can also be written as the braided tensor product $H_R\und\tens B$. (In the dual form in \cite{Ma:bos} the coproduct actually arises as braided tensor product with the adjoint action on $H$ and is then converted to bosonic form by an inverse transmutation). We are therefore in the setting of (\ref{tensbund}). \endproof

We see that every morphism $A:B\to \Omega^1H$ with $A(1)=0$ provides a connection by
\eqn{tensconn}{ \omega(b)=1\tens \und S b\Bo \tens 1\tens b\Bt -1\tens 1\tens 1\tens 1\eps(b)+ \Psi(\und S b\Bo\tens A(b\Bt)_i) \tens A(b\Bt)^i \tens b\Bt}
where $A(b)=A(b)_i\tens A(b)^i\in \Omega^1H\subset H\tens H$ and $\Psi(b\tens h)=h\o\tens b\bo\CR(b\bt\tens h\t)$ is the braiding between $B$ and $H_R$. This connection provides, in turn, a splitting of $\Omega^1$ of $H\rbiprod B$, as well as a covariant derivative on associated braided fiber bundles. We now underline the braided coproduct and antipode to distinguish them from usual Hopf algebra structures.

On the other hand, we can also view any bosonisation, being a semidirect product algebra, as a quantum principal bundle, as explained (with emphasis on the quantum double) in \cite{BrzMa:gau}. From this point of view the gauge quantum group is $H\subset H\rbiprod B$ and $B$ is the base. In this case every map $A:H\to \Omega^1B$ with $A(1)=0$ provides a connection by
\eqn{qgtensconn}{ \omega(h)=Sh\o\tens 1\tens h\t\tens 1-1\tens 1\tens 1\tens 1\eps(h)+Sh\o\tens A(h\t)_i\tens h\th\tens A(h\t)^i}
where $A(h)=A(h)_i\tens A(h)^i\in \Omega^1B$. We see that bosonisations are bundles in two ways:
\ceqn{twobundles}{  H\rbiprod B\\
\nearrow  \qquad\quad\qquad \nwarrow\\
 H \qquad\qquad\qquad\qquad  B}
and that gauge fields in one case are a certain class of maps $B\to H\tens H$ and in the other case a certain class of maps $H\to B\tens B$.

To give a concrete example of interest in physics, let $H=\widetilde{SO_q(n)}=\<\vect,\dila\>$ the dilaton-extended rotation group of any type associated to an $R$-matrix $R\in M_n\tens M_n$. Let $B=\R_q^n=\<p_i\>$ the braided plane associated to a compatible matrix $R'\in M_n\tens M_n$ \cite{Ma:poi}. Let $P=\widetilde{SO_q(n)}\lbiprod \R_q^n$ be generated by these as in \cite{Ma:poi}, to form the extended $q$-Poincar\'e function algebra as introduced there. Then we can view this as a braided $\R_q^n$-bundle over $\widetilde{SO_q(n)}$. A braided group gauge field is a function on $\R_q^n$ with values in $\Omega^1 \widetilde{SO_q(n)}$, i.e. provided by evaluation against an element of $\R_q^n\tens\Omega^1\widetilde{SO_q(n)}$ since $\R_q^n$ is self-dual via the quantum metric. Although we have given the details here only for the universal calculus, one \cite{BrzMa:coa} tells us also how to make quotients of the above constructions (details for the diagrammatic-braided case will be given elsewhere) and thereby work with non-universal calculi as well. For realistic calculi close to the classical ones, we can think of the vector space of $\Omega^1(\widetilde{SO_q(n)})$  as something like $SO_q(n)$ tensor a vector space of generators of $U_q(so_n)$. 

Moreover, we can also view this same $q$-algebra in a more conventional way as a trivial quantum group principal bundle over $\R_q^n$ with structure quantum group $\widetilde{SO_q(n)}$. Here the trivialisation is the Hopf algebra inclusion $\widetilde{SO_q(n)}\subset \widetilde{SO_q(n)}\rbiprod \R_q^n$. A quantum group gauge field in this case is a function on $\widetilde{SO_q(n)}$ with values in $\Omega^1\R_q^n$, i.e. an element of $\widetilde{U_q(so_n)}\tens \Omega^1\R_q^n$. In the non-universal case $\Omega^1(\R_q^n)$ as a vector space can be identified with $\R_q^n$ tensor a vector space of generators of the fermionic quantum plane $\R_q^{|n}$. 
 
Hence we have an example of a double fibration
\ceqn{twobundex}{ \widetilde{SO_q(n)}\rbiprod \R_q^n\\
\nearrow \qquad\qqquad\quad\qquad \nwarrow\\
\widetilde{SO_q(n)} \qquad\qqquad\qquad\qquad\qquad \R_q^n}
and a gauge field in the two cases actually involves rather similar data. This
suggests an interesting direction for further research.

Finally, we add a third point of view, as an embeddable homogeneous space bundle in the setting of the more general coalgebra gauge theory or $\psi$-bundles  in \cite{BrzMa:coa}. The data for an emebeddable homogeneous space bundle is a Hopf algebra $\CH$, a coalgebra $C$ and a coalgebra surjection  $\pi:\CH\to C$ such that $\ker\pi$ is a  right ideal subject to a minimality requirement. One may then construct an {\em entwining structure} $\psi:C\tens \CH\to \CH\tens C$ by $\psi(c\tens h)=h\o\tens \Pi(g h\o)$ for $g\in\pi^{-1}(c)$. The entwining structure induces a coaction $\cora(h)=\psi(\pi(1)\tens h)$ and under a suitable minimality requirement we have the corresponding $\chi$ (defined as before) invertible. See \cite{BrzMa:gau} for details.  

\begin{propos} Every bosonisation in the braided category of right $H$-comodules can be viewed as  an emebeddable homogeneous space bundle over the coalgebra $B$ with $\pi:H\rbiprod B\to B$ defined by $\pi(h\tens b)=\eps(h)b$. 
\end{propos}
\proof We set $\CH=H\rbiprod B$, the bosonisation, and $C=B$. Because the coproduct is a semidirect one by the right coaction of $H$ on $B$, it is easy to see that $\pi(h\tens b)=\eps(h)b$ is a coalgebra map.  Because the product is a semidirect one by an induced right action of $H$, $\ra$ say, we have $\pi((g\tens c)(h \tens b))=\pi(g h\o\tens (c\ra h\t)b)=\eps(g)(c\ra h)b=(\pi(g\tens c)\ra h) b$ so that $\ker\pi$ is a right ideal. Then $\psi$ has the form 
\align{\psi(c\tens(h\tens b))\equad &&=h\o\tens b\Bo\bo\tens \pi((1\tens c)(h\t b\Bo\bt\tens b\Bt))\\
&&=h\o\tens b\Bo\bo\tens (c\ra (h\t b\Bo\bt))b\Bt=h\o\tens b\Bo\bo\tens c\bo b\Bt\CR(c\bt\tens h\t b\Bo\bt)}
where we choose representative $1\tens c\in\pi^{-1}(c)$ and put in the explicit form of $\ra$ from bosonisation theory. The induced coaction is $\psi(1\tens (h\tens b))=h\tens b\Bo\tens b\Bt$. With more work, one may prove the invertibility of the induced $\chi$. However, this action and $\psi$ coincide with the $\psi$-bundle corresponding to the braided principal bundle in Proposition~4.1, for which the $\chi$ condition is known, so we do not need to
prove this. To see that the $\psi$ bundles are the same we need the formulae in \cite{BrzMa:coa} which construct a $\psi$ bundle from a braided principal bundle, namely $\psi(c\tens u)=\Psi(c\tens u\bo)u\bt$ for $c\in B$ and $u\in P$. In our case, the $B$-coaction is $\id\tens\Delta$ so 
$\psi(c\tens (h\tens b))=\Psi(c\tens (h\tens b\Bo))b\Bt=(h\tens b\Bo)\bo\tens c\bo b\Bt \CR(c\bt\tens (h\tens b\Bo)\bt)$ where $(h\tens b)\bo\tens (h\tens b)\bt=h\o\tens b\bo\tens h\t b\bt$ is the coaction of $H$ on $H_R\tens B$ as an object in our braided category of comodules (which has braiding $\Psi$ induced via the coactions from $\CR$). Putting this in, we obtain the same $\psi$ as for the embeddable homogeneous space bundle.  \endproof

\baselineskip 13pt  

\end{document}